\documentclass[letterpaper, 10 pt, conference]{ieeeconf}  

\IEEEoverridecommandlockouts                              

\usepackage{graphics} 
\usepackage{epsfig} 
\usepackage{times} 
\usepackage{amsthm}
\usepackage{mathtools}
\usepackage{amsmath} 
\usepackage{amssymb}  
\usepackage{color}

\usepackage{tikz}
\usetikzlibrary{shapes,arrows}
\usepackage{pgfplots} 
\usepackage{pgfgantt}
\usepackage{pdflscape}
\usepackage{nicefrac}
\usepackage{prettyref}
\pgfplotsset{compat=1.5}
\pgfplotsset{plot coordinates/math parser=false}
\newlength\fwidth
\usepackage{grffile}
\usetikzlibrary{plotmarks}
\usetikzlibrary{arrows.meta}
\usepgfplotslibrary{patchplots}
\usepackage{units}

\newtheorem{definition}{Definition}
\newtheorem{remark}{Remark}
\newtheorem{theorem}{Theorem}

\newtheorem{problem}{Problem}
\newtheorem{proposition}{Proposition}
\newtheorem{assumption}{Assumption}

\newcommand{\defeq}{\vcentcolon=}
\newcommand{\T}{\scriptscriptstyle\top}       

\newcommand{\mathmax}{\operatorname*{max}}
\newcommand{\mathargmin}{\operatorname*{arg\,min}}
\newcommand{\mathst}{\text{s.t.}}

\definecolor{myred}{RGB}{233,72,73}%
\definecolor{mygreen}{RGB}{113,191,110}%
\definecolor{myblue}{RGB}{93,147,191}%
\definecolor{mydarkblue}{RGB}{57,101,181}%
\definecolor{mycyan}{rgb}{ 0.05, 0.80, 0.75}%

\newrefformat{sec}{\mbox{Section \ref{#1}}}
\newrefformat{tab}{\mbox{Tab. \ref{#1}}}
\newrefformat{fig}{\mbox{Fig. \ref{#1}}}
\newrefformat{rem}{\mbox{Remark \ref{#1}}}
\newrefformat{alg}{\mbox{Algorithm \ref{#1}}}
\newrefformat{ass}{\mbox{Assumption \ref{#1}}}
\newrefformat{def}{\mbox{Definition \ref{#1}}}
\newrefformat{prop}{\mbox{Proposition \ref{#1}}}
\newrefformat{prob}{\mbox{Problem \ref{#1}}}
\newrefformat{thm}{\mbox{Theorem \ref{#1}}}

\renewenvironment{proof}{{\bfseries Proof:~}}{}

 \title{\LARGE \bf
 	Control-sharing Control Barrier Functions for Intersection Automation under Input Constraints}

\author{Alexander Katriniok
\thanks{Alexander Katriniok is with the Ford Research \& Innovation Center (RIC), 
	S\"{u}sterfeldstr. 200, 52072 Aachen, Germany, {\tt\small de.alexander.katriniok@ieee.org}}%
}

\begin{document}

\maketitle
\thispagestyle{empty}
\pagestyle{empty}

\begin{abstract}
This contribution introduces a centralized input constrained optimal control framework based on multiple control barrier functions (CBFs) to coordinate connected and automated agents at intersections. For collision avoidance, we propose a novel CBF which is safe by construction. The given control scheme provides provable guarantees that collision avoidance CBFs and CBFs to constrain the agents' velocity are jointly feasible (referred to as control-sharing property) 
subject to input constraints. A simulation study finally provides evidence that the proposed control scheme is safe.  
\end{abstract}

\section{Introduction}

Automated vehicles with communication capabilities are an 
emerging trend in automotive industry which gives rise to many challenging control problems. One of those problems, that has very actively been researched within last two decades, is the automation of road intersections. 
There is a rich body of literature on methods that have been applied to tackle the control problem at hand \cite{Khayatian2020a}. 
Optimization-based methods like model predictive control are an appealing choice for such kind of motion planning problems to explicitly accommodate constraints or exploit anticipated trajectories of conflicting agents \cite{Malikopoulos2021a,Katriniok2019a}.  
However, it is oftentimes challenging to solve the underlying optimization problems in real-time, especially when they are nonconvex.
An attractive alternative and computationally efficient way to handle constraints in an optimal control framework are control barrier functions (CBFs), which have to satisfy Lyapunov-like conditions to render a safe set forward invariant \cite{Ames2016a,Ames2019a}. 

While our previous research work has primarily built upon model predictive control (MPC) \cite{Katriniok2019a,Katriniok2017a}, this contribution aims to evaluate the potential of utilizing CBFs for intersection automation. Given that CBFs have also their limitations and drawbacks \cite{Ames2019a,Wang2017a}, our overall motivation (which goes far beyond the scope of this paper) is to investigate the pros and cons of both concepts on the example of intersection automation and to potentially come up with a suitable combination of the two.
To reasonably limit our scope, we start with a centralized control approach where a central node within the infrastructure determines agent controls and ensures collision avoidance by solving a quadratic programming (QP) problem subject to linear CBF constraints. The resulting control actions are then sent back to the agents via vehicle-to-infrastructure (V2I) communication. 

In literature, CBFs have frequently been used to avoid collisions between aircrafts \cite{Squires2018a} or vehicles \cite{Ames2016a,Ames2019a,Xiao2019a,Santillo2021a}. Intersection automation using CBFs, though, has not so frequently been investigated in the community, see \cite{Khaled2020a,Shivam2020a} for some recent contributions. The utilization of CBFs, though, gets challenging when input constraints are present \cite{Ames2019a,Agrawal2021a} or multiple CBFs \cite{Xu2018a} are applied simultaneously. Then, it may happen that the resulting QP gets infeasible. To maintain safety guarantees, though, the CBFs have to be designed such that they are jointly feasible subject to input constraints \cite{Ames2019a,Squires2018a,Agrawal2021a}.

\subsection{Contribution}

In the remainder, we convey a centralized input constrained optimal control framework which relies on multiple CBFs 
to satisfy state and safety constraints, subject to bounded control inputs. Our contribution extends the state-of-the-art as follows. 
Firstly, we propose a novel CBF for collision avoidance that is safe by construction. To this end, we build upon the general idea in \cite{Ames2016a} and extend it towards the 2-dimensional case of a superellipse. Secondly, we utilize multiple CBFs simultaneously, that is, to ensure collision avoidance and to constrain the agents' velocity. For such setup, we prove that these CBFs are jointly feasible --- referred to as control-sharing property \cite{Xu2018a}. Thirdly, the utilized CBFs are provably compatible with input constraints. The latter two contributions actually constitute two of the major challenging aspects of applying CBFs for constrained control. Lastly, we contribute with novel ideas and perspectives to the rather sparse literature on CBF-based intersection automation.

While the approach in \cite{Shivam2020a} is entirely different from ours, \cite{Khaled2020a} also builds upon the idea in \cite{Ames2016a} to design a CBF that is safe by construction. However, we propose a less conservative formulation for collision avoidance, we contemplate the actual relative motion of the agents in the horizontal plane and we provide provable guarantees for the CBFs to be control-sharing subject to input constraints.

In the remainder, \prettyref{sec:problem} outlines the intersection coordination problem. Then, \prettyref{sec:CBF} and \prettyref{sec:CBFCA} convey the construction of CBFs and show that these are provably safe and compatible with input constraints. Thereafter, controller synthesis is introduced in \prettyref{sec:synthesis} before simulation results are discussed in \prettyref{sec:sim}. 

\section{Intersection Coordination Problem}
\label{sec:problem}

\subsection{Problem Description}
\label{sec:problem_description}

The intersection automation problem to solve involves a set \mbox{$\mathcal{A} \defeq \{1,\ldots,N_\mathcal{A}\}$} of connected and automated agents, and can be stated as follows.
\begin{problem} \normalfont \label{prob:problem_description_problemStatement}
	The objective is to automate agents at a four way, single lane intersection without any traffic lights or signs by manipulating their acceleration. 
	A central node in the infrastructure determines the control actions of agents \mbox{$i \in \mathcal{A}$} such that they track their reference speed and cross the intersection without colliding with other agents $j \in \mathcal{A}\setminus\{i\}$.
\end{problem}
\vspace*{-0.5mm}
To limit the complexity of the control problem at hand, we make the following assumptions.
\vspace*{-0.5mm}
\begin{assumption}	\normalfont
	\mbox{A1. All agents} cross the intersection straight; 
	\mbox{A2. Rear-end} collisions are out of scope;
	A3. Perfect V2I communication is used for information exchange (no packet losses, no communication delays);
	\mbox{A4. Agent} models are not subject to uncertainties and agent states are fully measurable. 
	\mbox{A5.} Lateral vehicle control is accommodated by a separate and independent control layer.
	\label{ass:problem_description_problemAssump}
\end{assumption}

\begin{figure}[b!]
	\centering
	\def\svgwidth{60mm}
	\vspace*{-4mm}
	\input{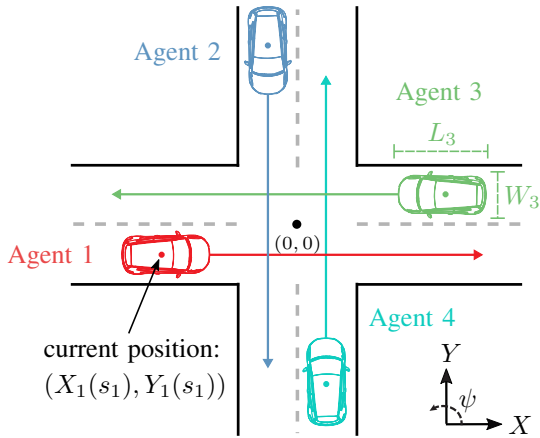}	
	\vspace*{-4mm}
	\caption{Scenario with 4 agents crossing the intersection straight.}%
	\label{fig:problem_description_intersectionScenario}
\end{figure}

\subsection{Modeling} 
\label{sec:problem_modeling}

We devise the dynamics of every agent $i \in \mathcal{A}$ as a nonlinear input-affine state space model of the form
\begin{align}
	\underbrace{\vphantom{
			\begin{bmatrix}
			-F_{r,i}(v_i)/m \\
			v_i		
			\end{bmatrix}}
	\begin{bmatrix}
	\dot{v_i} \\
	\dot{s_i}
	\end{bmatrix}}_{\dot{x}_i} = 
	\underbrace{\begin{bmatrix}
	  -F_{r,i}(v_i)/m_i \\
	  v_i
	\end{bmatrix}}_{f_i(x_i)} + 
	\underbrace{\vphantom{
			\begin{bmatrix}
			-F_{r,i}(v_i)/m \\
			v_i		
			\end{bmatrix}}
   \begin{bmatrix}
	1 \\
	0
	\end{bmatrix}}_{g_i(x_i)} 	
    \underbrace{\vphantom{
      \begin{bmatrix}
	  -F_{r,i}(v_i)/m \\
      v_i		
      \end{bmatrix}} a_{x,i} }_{u_i}
  \label{eq:problem_modeling_xdot}
\end{align}
which describes the evolution of Agent $i$'s velocity $v_i$ and path coordinate $s_i$, influenced by the manipulated longitudinal acceleration $u_i=a_{x,i}$, along  an \textit{a priori} given path $\mathcal{P}_i$. The resistance force $F_{r,i}(v_i)$ in \eqref{eq:problem_modeling_xdot}, which involves rolling resistances and aerodynamic drag, is approximated as a quadratic function of the agent's velocity \cite{Rajamani2012}, i.e., 
\begin{align*}
F_{r,i}(v_i) = \mathrm{sign}(v_i) c_{0,i} + c_{1,i} v_i + c_{2,i} v_i^2
\end{align*}
with $c_{0,i},\,c_{1,i},\,c_{2,i} \in \mathbb{R}$. 
Moreover, $m_i \in \mathbb{R}_{>0}$ is the mass of \mbox{Agent $i$}. For notational convenience, we also stack the state space models \eqref{eq:problem_modeling_xdot} of agents $i$ with state vector $x_i = [v_i,s_i]^{\T}$ and control input $u_i = a_{x,i}$, such as to come up with an augmented state space model of the form
\begin{align}
\dot{x} = f(x) + g(x)u
  \label{eq:problem_modeling_lumpedSS}
\end{align}
where $x=\{x_i\}_{i=1}^{N_\mathcal{A}}$, $u=\{u_i\}_{i=1}^{N_\mathcal{A}}$ are the stacked state and input vectors, and $f=\{f_i\}_{i=1}^{N_\mathcal{A}}$ and $g=\{g_i\}_{i=1}^{N_\mathcal{A}}$ are the stacked vector fields respectively.

\subsection{Local and Global Coordinates}
\label{sec:problem_coordinates}
To infer the global Cartesian coordinates $(X_i, Y_i)$ of every agent $i$ from its local path coordinate $s_i$, we utilize the mapping function \mbox{$\mathcal{P}_i: s_i \mapsto (X_i, Y_i)$}.
For the considered straight crossing use case in \prettyref{fig:problem_description_intersectionScenario}, it suffices to utilize linear affine mapping functions, that is,
\begin{align*}
	X_i(s_i) \defeq p_{x,i,0} + p_{x,i,1} s_i, ~~ Y_i(s_i) \defeq p_{y,i,0} + p_{y,i,1} s_i. 
\end{align*}
with constants $p_{x,i,0}$, $p_{x,i,1}$, $p_{y,i,0}$, $p_{y,i,1} \in \mathbb{R}$.
 Depending on the use case, these mapping functions can of course be more complex, e.g., by applying B-Spline functions 
 to represent the paths of turning agents \cite{Katriniok2019a}. Likewise, every agent's yaw angle can be described as a function $\psi_i(s_i)$ of its path coordinate $s_i$.  With \prettyref{ass:problem_description_problemAssump}, for the straight crossing use case holds $\psi_i(s_i) = n\frac{\pi}{2}$ with $n=[0,3] \subset \mathbb{Z}$.


\section{CBFs for Constrained Control}
\label{sec:CBF}

\subsection{Background on CBFs}
\label{sec:CBF_backgrnd}

We consider input-affine dynamical systems whose state evolution can be described according to \eqref{eq:problem_modeling_lumpedSS}
with locally Lipschitz vector fields $f$ and $g$, state vector $x \in 
\mathbb{R}^{n_x}$ and vector of control inputs $u \in 
\mathbb{R}^{n_u}$.
Along these lines, we define the zero superlevel set of a continuously differentiable function \mbox{$h: \mathcal{D} \rightarrow \mathbb{R}$} as
\begin{align}
\mathcal{C} &\defeq \left\{ x \in \mathcal{D} \mid h(x) \geq 0 \right\}, \notag \\
\delta\mathcal{C} &\defeq \left\{ x \in \mathcal{D} \mid h(x) = 0 \right\}, \label{eq:CBF_backgrnd_Cdef} \\
\mathrm{Int}(\mathcal{C}) &\defeq \left\{ x \in \mathcal{D} \mid h(x) > 0 \right\}. \notag 
\end{align}
\begin{definition}[\cite{Ames2019a}]
	The set $\mathcal{C}$ is said to be forwards invariant if $x(t_0) \in \mathcal{C}$ implies $x(t) \in \mathcal{C}$ for all $t\geq t_0$.	
\end{definition}
\begin{definition}[\cite{Ames2016a}] Let $\mathcal{C} \subseteq \mathcal{D} \subset \mathbb{R}^{n_x}$ be a zero superlevel set of a continuously differentiable function  $h: \mathcal{D} \rightarrow \mathbb{R}$ in accordance to \eqref{eq:CBF_backgrnd_Cdef}. For a dynamical system of the form \eqref{eq:problem_modeling_lumpedSS}, $h$ is a control barrier function (CBF) if there exists an extended class $\mathcal{K}_\infty$ function \mbox{$\alpha: \mathbb{R} \rightarrow \mathbb{R}$} such that the inequality 
	\vspace*{-1mm}
	\begin{align}
	\sup_{u\in\mathcal{U}} [ L_f h(x) + L_g h(x) u + \alpha(h(x)) ] \geq 0 
	\label{eq:CBF_backgrnd_CBFdef} \\[-7mm] \notag
	\end{align}
	holds for all $x \in \mathcal{D}$. 
	In this context, $L_f h(x)= \frac{\partial h(x)}{\partial x}f(x)$ and $L_g h(x)= \frac{\partial h(x)}{\partial x}g(x)$ are the Lie derivatives of $h(x)$ along $f$ and $g$ respectively. 
	As \mbox{$h(x)=0$} for all $x \in \partial\mathcal{C}$, $h$ is also referred to as zeroing CBF (ZCBF).
	\label{def:CBF_backgrnd_CBF}
\end{definition}
In other words, $h$ is a CBF if there exists at least one admissible $u \in \mathcal{U}$ for every $x \in \mathcal{D}$ which renders \eqref{eq:CBF_backgrnd_CBFdef} feasible.

\begin{theorem}[\cite{Ames2019a}]
Let $\mathcal{C}\subseteq \mathbb{R}^{n_x}$ be a zero superlevel set 
of a continuously differentiable function \mbox{$h: \mathcal{D} \rightarrow \mathbb{R}$}. If $h$ is a CBF on $\mathcal{D} \supseteq \mathcal{C}$, then any Lipschitz continuous controller $u(x) \in K_{\text{CBF}}$ for the dynamical system \eqref{eq:problem_modeling_lumpedSS} with
\vspace*{-1mm}
\begin{align*}
	 K_{\text{CBF}}(x) \defeq \left\{ u \in \mathcal{U} : L_f h(x) + L_g h(x) u + \alpha(h(x)) \geq 0 \right\} \\[-6mm]\notag
\end{align*}   
renders the set $\mathcal{C}$ forward invariant and as such safe.
\label{thm:CBF_backgrnd_KCBF}
\end{theorem}
\begin{definition}[\cite{Xu2018a}]
 The CBF $h(x)$ for the dynamical system \eqref{eq:problem_modeling_lumpedSS} is said to have relative degree $r \in \mathbb{N}_{>0}$ if the $r$-th time derivative of $h(x)$ depends on the input $u$ while $1\leq r^\prime < r$ does not, that is, $L_g^r h(x) \neq 0$ while $L_g^{r^\prime} h(x) = 0$.
\end{definition}
For a CBF $h(x)$ with relative degree one and a linear extended class $\mathcal{K}_\infty$ function $\alpha: h(x) \mapsto \lambda h(x)$ with \mbox{$\lambda \in \mathbb{R}_{>0}$}, we can devise the differential inequality
\vspace*{-1mm}
\begin{align*}
\underbrace{L_f h(x) + L_g h(x) u}_{\dot{h}(x,u)} \geq -\lambda h(x). 
\end{align*}
By the Comparison Lemma \cite[Lemma 3.4]{Khalil2002a}, we can infer
\vspace*{-0.5mm}
\begin{align*}
h(t) \geq h(t_0) e^{-\lambda t}.\\[-5.5mm]
\end{align*}
In other words, if the system is safe at time $t_0=0$ with $h(t_0)\geq 0$, then it is also safe at time $t\geq t_0$, i.e., $h(t)\geq 0$. When applying multiple CBFs, say $h_i(x)$ with $i=1,\ldots,N_h$, it needs to be ensured that these are jointly feasible which is defined as control-sharing property \cite{Xu2018a}.
\begin{definition}
	The CBFs $h_i(x)$, $i=1,\ldots,N_h$ are said to have control-sharing property if the intersection of their set of controls is non-empty, that is, $\cap_{i=1}^{N_h} K_i(x)\neq \emptyset$.
\end{definition} 

\vspace*{1mm}
In the scope of intersection automation, we aim to \mbox{i) constrain} the agents' velocity to enforce forward driving with a maximum speed (see \prettyref{sec:CBF_velCBF}) and ii) to avoid collisions between agents (see \prettyref{sec:CBFCA}). Both conditions shall be imposed as CBF constraints of relative degree one, subject to bounds on the inputs, i.e.,
\begin{align}
u_i \in \mathcal{U}_i \defeq [ \underline{u}_i, \overline{u}_i ] =  [ \underline{a}_{x,i}, \overline{a}_{x,i} ] 
\label{eq:CBF_backgrnd_inputCons}
\end{align}
with $\underline{a}_{x,i} \in \mathbb{R}_{< 0}$ and $\overline{a}_{x,i} \in \mathbb{R}_{> 0}$. For such input constrained systems, it is a common approach to consider $\mathcal{D}$ in \prettyref{def:CBF_backgrnd_CBF} to be equal to $\mathcal{C}$ \cite{Squires2018a,Agrawal2021a}. In other words, the initial state $x(t_0)$ needs to be located within the safe set $\mathcal{C}$.


\subsection{CBFs for Constrained Agent Velocity}
\label{sec:CBF_velCBF}

The velocity of every agent $i$ shall be bounded between a minimum velocity $\underline{v}_i=0$ (only driving forward) and a maximum velocity $\overline{v}_i$, that is,
\begin{align}
\underline{v}_i \leq v_i \leq \overline{v}_i.
\label{eq:CBF_velCBF_upperLowerBnd}
\end{align}
These bounds can be formalized by the following two continuous differentiable functions
\vspace*{-1mm}
\begin{subequations}
	\begin{align}
	h_{\underline{v},i}(x_i) &\defeq ~~v_i - \underline{v}_i \label{eq:CBF_velCBF_hLowerBnd}\\  
	h_{\overline{v},i}(x_i) &\defeq -v_i + \overline{v}_i  \label{eq:CBF_velCBF_hUpperBnd}
	\end{align}
\end{subequations}
which need to be greater than or equal to zero to satisfy \eqref{eq:CBF_velCBF_upperLowerBnd}. We refer to the respective zero superlevel sets as $\mathcal{C}_{{\underline{v},i}}$ and $\mathcal{C}_{{\overline{v},i}}$.
By utilizing a linear extended class $\mathcal{K}_\infty$ function, the resulting CBF constraints (see \prettyref{thm:CBF_backgrnd_KCBF}) can be written as
\vspace*{-4mm}
\begin{subequations}
\begin{align}
	L_{f_i} {h}_{\underline{v},i}(x_i) + L_{g_i} {h}_{\underline{v},i}(x_i) u_i  + \lambda_{\underline{v}} {h}_{\underline{v},i}(x_i) &\geq 0 \label{eq:CBF_velCBF_CBFconsLowerBnd}  \\
	L_{f_i} {h}_{\overline{v},i}(x_i) + L_{g_i} {h}_{\overline{v},i}(x_i) u_i  + \lambda_{\overline{v}} {h}_{\overline{v},i}(x_i) &\geq 0
	\label{eq:CBF_velCBF_CBFconsUpperBnd}	
\end{align}
\end{subequations}
with positive constants $\lambda_{\overline{v}},\lambda_{\underline{v}} \in \mathbb{R}_{>0}$. We finally need to prove that \eqref{eq:CBF_velCBF_hLowerBnd} and \eqref{eq:CBF_velCBF_hUpperBnd} are actually CBFs on $\mathcal{C}_{{\underline{v},i}}$ and $\mathcal{C}_{{\overline{v},i}}$.
\begin{proposition}
		\label{prop:CBF_velCBF_propCBFvel}	
	The continuous differentiable functions \eqref{eq:CBF_velCBF_hLowerBnd} and \eqref{eq:CBF_velCBF_hUpperBnd} are CBFs on $\mathcal{C}_{{\underline{v},i}}$ and $\mathcal{C}_{{\overline{v},i}}$ for the dynamical system \eqref{eq:problem_modeling_lumpedSS}. That said, the sets $K_{{\underline{v},i}}(x_i)$ and $K_{{\overline{v},i}}(x_i)$ of controls s.t. $u_i\in \mathcal{U}_i$ (as defined in (5)), that render the zero superlevel sets $\mathcal{C}_{{\underline{v},i}}$ and $\mathcal{C}_{{\overline{v},i}}$ 
	forward invariant, are non-empty.
\end{proposition}
\begin{proof} We first contemplate the \textit{lower velocity bound}. By substituting the Lie derivatives in CBF constraint \eqref{eq:CBF_velCBF_CBFconsLowerBnd}, it holds 
	$a_{x,i}  \geq F_{r,i}(v_i)/m_i + \lambda_{\underline{v}} (\underline{v}_i-v_i)$.
It can be seen that CBF \eqref{eq:CBF_velCBF_hLowerBnd} (same for \eqref{eq:CBF_velCBF_hUpperBnd}) has relative degree one. With $\underline{v}_i=0$ and a technically reasonable parameterization of $F_{r,i}(v_i)$ and a suitable choice of $\overline{a}_{x,i}$, we can ensure  \mbox{$\overline{a}_{x,i} \geq \mathmax_{v_i}\{ F_{r,i}(v_i)/m_i - \lambda_{\underline{v}} {v}_i\}$}. That way, it holds
\begin{align}
\underbrace{\overline{a}_{x,i}}_{>0} \geq a_{x,i}  \geq \underbrace{F_{r,i}(v_i)/m_i}_{\geq 0} - \underbrace{\lambda_{\underline{v}} v_i}_{\geq 0}. 
\label{eq:CBF_velCBF_vminLowerAccel}
\end{align}
Thus, $K_{{\underline{v},i}}(x_i)$ is non-empty and $h_{\underline{v},i}$ is a CBF on $\mathcal{C}_{{\underline{v},i}}$. 
\hfill $\square$
\end{proof}

For the \textit{upper velocity bound} the proof follows the same reasoning, i.e., substituting the Lie derivatives in \eqref{eq:CBF_velCBF_CBFconsUpperBnd} gives 
\begin{align}
\underbrace{\underline{a}_{x,i}}_{<0} \leq a_{x,i}  \leq \underbrace{F_{r,i}(v_i)/m_i}_{\geq 0} + \underbrace{\lambda_{\overline{v}}(\overline{v}_i - v_i)}_{\geq 0}.
\label{eq:CBF_velCBF_vmaxUpperAccel}
\end{align}
Thus, $K_{{\overline{v},i}}(x_i)$ is non-empty and $h_{\overline{v},i}$ is a CBF on $\mathcal{C}_{{\overline{v},i}}$. \hfill $\square$


\subsection{Control-sharing Property}
We finally need to answer the question whether \eqref{eq:CBF_velCBF_hLowerBnd} and \eqref{eq:CBF_velCBF_hUpperBnd} are jointly feasible.
\begin{proposition}
	CBFs \eqref{eq:CBF_velCBF_hUpperBnd} and \eqref{eq:CBF_velCBF_hLowerBnd} exhibit a control-sharing property, that is, $K_{{\underline{v},i}}(x_i) \cap K_{{\overline{v},i}}(x_i) \neq \emptyset$ iff the inequality $-\lambda_{\underline{v}} v_i \leq -\lambda_{\overline{v}} v_i + \lambda_{\overline{v}} \overline{v}_i$ holds. The inequality can, e.g., be satisfied by choosing  $\lambda_{\underline{v}}\geq\lambda_{\overline{v}}$.
	\label{prop:CBF_velCBF_csProperty_prop}
\end{proposition}
\vspace*{-1mm}
\begin{proof}
	 If $-\lambda_{\underline{v}} v_i \leq -\lambda_{\overline{v}} v_i + \lambda_{\overline{v}} \overline{v}_i$  holds, \eqref{eq:CBF_velCBF_vminLowerAccel}-\eqref{eq:CBF_velCBF_vmaxUpperAccel}  imply
	\begin{align*}
	\hspace*{-2mm}
	\underbrace{F_{r,i}(v_i)/m_i}_{\geq 0} - \underbrace{\lambda_{\underline{v}} v_i}_{\geq 0} \leq a_{x,i}  \leq \underbrace{F_{r,i}(v_i)/m_i}_{\geq 0} - \underbrace{\lambda_{\overline{v}} v_i}_{\geq 0} + \underbrace{\lambda_{\overline{v}} \overline{v}_i}_{\geq 0}.
	\end{align*}
	Thus, we can infer that $K_{{\underline{v},i}}(x_i) \cap K_{{\overline{v},i}}(x_i) \neq \emptyset$ holds. \hfill$\square$
\end{proof}


\section{Collision Avoidance}
\label{sec:CBFCA}

\subsection{Geometric Formulation}
\label{sec:CBF_CA_geometry}

To avoid collisions between Agent $i$ and Agent $j$, we define a superellipse of forth degree, whose axes align with Agent $i$'s body frame $(X_i^b,Y_i^b)$ and whose origin coincides with \mbox{Agent $i$'s} geometric center \mbox{$P_i(s_i) \defeq (X_i(s_i),Y_i(s_i))$} (global frame), see 
\prettyref{fig:CBF_CA_formCBF_freebodyCBF}. That said, the superellipse is given by
\begin{align}
SE_{i}(P_i^b;a_{ij},b_{ij}) \defeq \frac{(X_i^b)^4}{a_{ij}^4} + \frac{(Y_i^b)^4}{b_{ij}^4} - 1
\label{eq:CBF_CA_geometry_SE}
\end{align}
where the axis length and width $a_{ij}, b_{ij} \in \mathbb{R}_{>0}$ are defined with respect to the width $W_i,W_j$ and length $L_i,L_j$ of agents $i$ and $j$, as well as their relative yaw angle $\psi_j-\psi_i$, i.e.,
\begin{align*}
	a_{ij} &\defeq \frac{L_i}{2} \,+ \| R(\psi_j-\psi_i)\begin{bmatrix}
		\,\nicefrac{L_j}{2}\, \\ 0
	\end{bmatrix} \|_2 + \Delta a_{ij} \\
	b_{ij} &\defeq \frac{W_i}{2} + \| R(\psi_j-\psi_i)\begin{bmatrix}
		0 \\ \nicefrac{W_j}{2}
	\end{bmatrix} \|_2 + \Delta b_{ij}.
\end{align*}
In that regard, the rotation matrix
\begin{align*} 
R(\psi) \defeq \begin{bmatrix}
	\cos(\psi) & -\sin(\psi) \\
	\sin(\psi) & ~\,\cos(\psi)
\end{bmatrix} 
\end{align*}
accommodates the relative yaw angle while $\Delta a_{ij},$ \mbox{$\Delta b_{ij} \in \mathbb{R}_{\geq 0}$} are additional safety buffers. Every point \mbox{$P_i^b= (X_i^b,Y_i^b)$} on the boundary of the superellipse satisfies $SE_{i}(P_i^b;a_{ij},b_{ij})=0$. That said, Agent $i$ and Agent $j$ are safe if the geometric center \mbox{$P_j(s_j) \defeq (X_j(s_j),Y_j(s_j))$} (global frame) of \mbox{Agent $j$} is either located outside or on the boundary of the superellipse, which is equivalent to 
\begin{align*}
SE_{ij}(x) \defeq SE_{i}\bigl(  R(\psi_i)^{\T}(P_j(s_j)-P_i(s_i);~a_{ij},b_{ij})  \bigr) \geq 0 
\end{align*}
where $R(\psi_i)$ transfers the vector $P_j(s_j)-P_i(s_i)$ into the body frame of Agent $i$. 
\subsection{CBF Formulation}
\label{sec:CBF_CA_formCBF}

\begin{figure}[b!]
	\centering
	\def\svgwidth{74mm}
	\vspace*{-4mm}
	\input{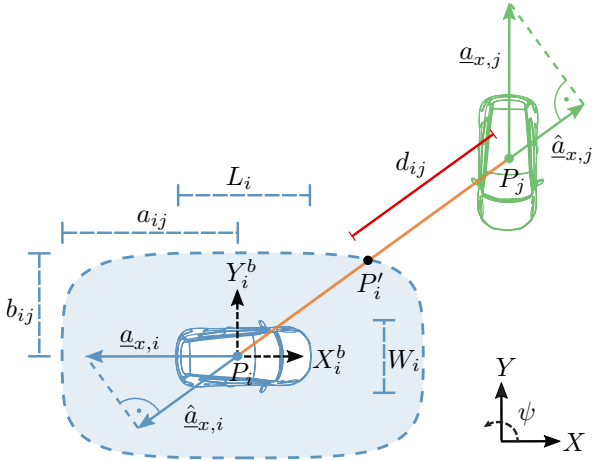}	
	\vspace*{-2mm}
	\caption{Freebody diagram illustrating the construction of the CBF in terms of i) the distance between Agent $i$'s superellipse and the geometric center of Agent $j$ and ii) the 
	maximum deceleration along the vector $P_j-P_i$.} 
	\label{fig:CBF_CA_formCBF_freebodyCBF}
\end{figure}
Along these lines, an intuitive approach would be to apply $h_{c,ij}(x)=SE_{ij}(x)$ as CBF.
Such selection of $h_{c,ij}(x)$ results in a CBF of relative degree 2 as the path coordinate needs to be differentiated twice until the control input $u$ appears in the CBF constraint. It is easy to prove, though, that $SE_{ij}(x)$ does not render the respective superlevel set $\mathcal{C}_{c,ij}$ forward invariant under input constraints $u \in \mathcal{U}$ with $\mathcal{U} \defeq \{ u\in\mathbb{R}^{N_\mathcal{A}} \mid \forall i: u_i \in \mathcal{U}_i \}$. More precisely, consider the case that $P_j$ in \prettyref{fig:CBF_CA_formCBF_freebodyCBF} is located on the boundary of the superellipse, i.e., $SE_{ij}(x)=0$. If the velocity of at least one agent is non-zero, a violation of safety constraints is unavoidable. If $u$ is unbounded, however, $SE_{ij}(x)$ may be a valid CBF as an arbitrarily large control input can be applied to avoid the collision. There are generally two ways to accommodate input constraints:
\begin{enumerate}
   \item By construction the CBF $h_{c,ij}(x)$ is forward invariant on $\mathcal{C}_{c,ij}$ \mbox{s.t. $u\in\mathcal{U}$} and as such safe \cite{Ames2019a,Breeden2021a}.
   \item Start with a CBF candidate and \textit{tighten} $\mathcal{C}_{c,ij}$ iteratively until the resulting CBF $h^\prime_{c,ij}(x)$ is forward invariant on $\mathcal{C}^\prime_{c,ij}(x)$ s.t. $u\in\mathcal{U}$ \cite{Ames2019a,Squires2018a,Agrawal2021a}.
\end{enumerate}
We pursue option 1 as it is more intuitive in its design. Along these lines, we build upon the general idea in \cite{Ames2016a} and extend it towards the 2-dimensional case of a superellipse. To this end, we define the CBF such that the distance $d_{ij}$ between the superellipse and \mbox{Agent $j$'s} geometric center 
is lower bounded by a safety distance $d_{\text{safe},ij}$ which allows both agents to stop without a collision 
when exploiting their maximum braking capability along the vector $P_j-P_i$, see \prettyref{fig:CBF_CA_formCBF_freebodyCBF}. 
\begin{remark}
	Such definition of the CBF does not imply that both agents have to stop or end up in a deadlock. However, CBFs are generally known to be prone to deadlock-like situations, which can however be circumvented \cite{Wang2017a}.
\end{remark}
To start with, we define 
\vspace*{-1mm}
\begin{align}
\hspace*{-2.5mm} P_i^\prime(s_i,s_j) \defeq P_i(s_i) + \nu \,\frac{P_j(s_j)-P_i(s_i)~}{ \| P_j(s_j) - P_i(s_i) \|_2},~\nu \in \mathbb{R}_{>0} \raisetag{3mm} \\[-7mm] \notag
\end{align}
as the point (global coordinates) on the boundary of the superellipse along the vector $P_j-P_i$, see \prettyref{fig:CBF_CA_formCBF_freebodyCBF}. By solving 
\begin{align*}
 SE_i\bigl(R(\psi_i)^{\T}(P_i^\prime(s_i,s_j)-P_i(s_i);~a_{ij},b_{ij}\bigr)=0
\end{align*}
for $\nu$, a positive real $\nu$ can uniquely be determined.
There are four roots, two of which are conjugate complex and two are real. Only one of the latter is positive and satisfies the condition $\nu \in \mathbb{R}_{>0}$. Finally, for the distance $d_{ij}$ between $P_j$ and the superellipse at $P_i^\prime$ holds
\vspace*{-1mm}
\begin{align}
d_{ij} \defeq \| P_j(s_j) - P_i(s_i) \|_2 - \| P_i^\prime(s_i,s_j) - P_i(s_i) \|_2 
\end{align}
\begin{remark}
The way we defined $d_{ij}$ 
ensures that $d_{ij}$ is negative when Agent $j$ is located inside the superellipse, as opposed to defining $d_{ij}$ as $\| P_j - P_i^\prime \|_2$.
\end{remark}
\begin{figure*}
	\centering
	\setlength\fwidth{0.60\textwidth}	
	\input{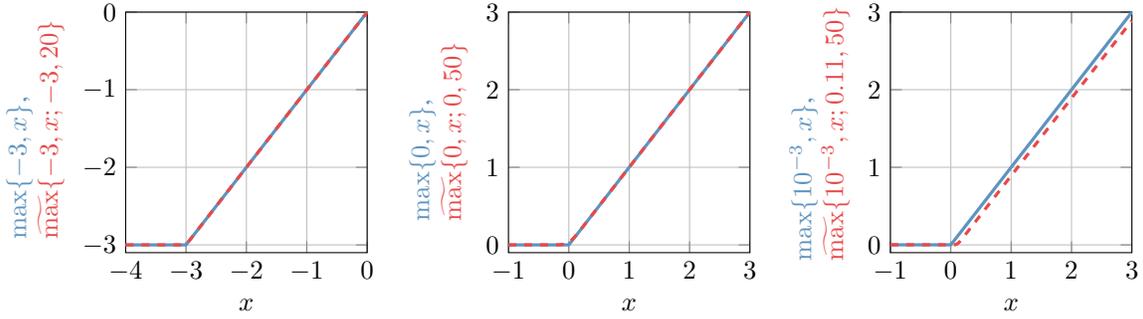}
	\vspace*{-4mm}
	\caption{Left: Smooth over-approximation ({\color{myred} dashed red}) of \eqref{eq:CBF_CA_formCBF_axMinEff} ({\color{myblue} solid blue}) for $\underline{a}_{x,i}=-3$, which under-approximates the resulting (positive) deceleration in \eqref{eq:CBF_CA_formCBF_axProj}; Middle: Smooth over-approximation ({\color{myred} dashed red}) in \eqref{eq:CBF_CA_formCBF_dSafe} ({\color{myblue} solid blue}); Right: Smooth under-approximation ({\color{myred} dashed red}) in \eqref{eq:CBF_CA_formCBF_axMinRelative} ({\color{myblue} solid blue}).}
	\vspace*{-6mm}
	\label{fig:problem_funcapprox}
\end{figure*}
The next step is to deduce a suitable safety distance $d_{\text{safe},ij}$ which renders $\mathcal{C}_{c,ij}$ forward invariant by construction. To that end, we introduce the time derivative of $d_{ij}$ 
\begin{align}
\hspace*{-1.7mm}{v}_{ij} \defeq \dot{d}_{ij}  = &\frac{(P_j - P_i)^{\T}}{\| P_j - P_i \|_2} \left( \frac{dP_j}{ds_j} v_j - \frac{dP_i}{ds_i} v_i \right) - \\
&\frac{(P_i^\prime - P_i)^{\T}}{\| P_i^\prime - P_i \|_2} \left( \frac{dP_i^\prime}{ds_i} v_i + \frac{dP_i^\prime}{ds_j} v_j - \frac{dP_i}{ds_i} v_i \right) \notag
\end{align}
which is negative when both agents are getting closer and positive when moving away from each other. Moreover, the maximum deceleration of Agent $i$ 
along $P_i - P_j$ (likewise with $P_j - P_i$ for \mbox{Agent $j$)} can be represented as 
\begin{align} \label{eq:CBF_CA_formCBF_axProj}
\underline{\hat a}_{x,i} \defeq \frac{(P_i - P_j)^{\T}}{ \| P_i - P_j \|_2} R(\psi_i) \begin{bmatrix}
\underline{a}_{x,\text{eff},i} \\ 0
\end{bmatrix} ~~ \\
\text{where} ~~~~~~\underline{a}_{x,\text{eff},i} \defeq \max\{ \underline{a}_{x,i},\,-\lambda_{\underline{v}} v_i \}
\label{eq:CBF_CA_formCBF_axMinEff} \\[-6mm] \notag
\end{align}
is the effective maximum deceleration (incl. resistances) of Agent $i$ which is compatible with the lower velocity bound \eqref{eq:CBF_velCBF_vminLowerAccel}. Without consideration of \eqref{eq:CBF_velCBF_vminLowerAccel}, the CBFs will not be control-sharing. 
By neglecting the resistance $-F_{r,i}(v_i)/m_i$ in the first term of \eqref{eq:CBF_CA_formCBF_axMinEff}, we are slightly more conservative. However, we can exploit the smooth approximation \eqref{eq:CBF_CA_formCBF_maxApprox} of the $\max$-operator, which requires one term to be constant.
The resulting safety distance between Agent $i$ \mbox{and Agent $j$} 
\begin{align}
d_{\text{safe},ij} \defeq \frac{\max\{0,-v_{ij}\}^2}{2\, \underline{\hat a}_{x,ij}} ~~~~~~~~~~~
\label{eq:CBF_CA_formCBF_dSafe} \\
\text{with} ~~~~~~ \underline{\hat a}_{x,ij} \defeq \max\{ \epsilon, \underline{\hat a}_{x,i} \} + \max\{ \epsilon, \underline{\hat a}_{x,j} \},
\label{eq:CBF_CA_formCBF_axMinRelative} \\[-6mm] \notag
\end{align}
and sufficiently small constant $\epsilon > 0$ 
allows us to define a collision avoidance CBF that is safe by construction (both agents are able to stop). The $\max$-operator in \eqref{eq:CBF_CA_formCBF_axMinRelative} ensures $\underline{\hat a}_{x,ij} > 0$ and avoids division by zero in \eqref{eq:CBF_CA_formCBF_dSafe}
when the agents are passing the intersecting point along their paths and are moving away from each other. The $\max$-operator in the nominator of \eqref{eq:CBF_CA_formCBF_dSafe}
causes $d_{\text{safe},ij}$ to be zero if the distance $d_{ij}$ between Agent $i$ and Agent $j$ is increasing. By using the $\max$-operator in \eqref{eq:CBF_CA_formCBF_axMinEff}--\eqref{eq:CBF_CA_formCBF_axMinRelative}, $d_{\text{safe},ij}$ and as such the resulting CBF is no longer a continuous differentiable function, as opposed to \prettyref{def:CBF_backgrnd_CBF}.
To this end, we introduce the smooth approximation of $\max\{ c, x \}$
\vspace*{-1mm}
\begin{align}
 \widetilde{\max}\{ c, x; b_1, b_2  \} &\defeq c+\frac{1}{b_2}\ln(1+e^{(x-b_1)b_2}) \label{eq:CBF_CA_formCBF_maxApprox} \\[-6mm] \notag  
\end{align}
with $c=\mathtt{const.}$, $b_1,b_2 \in \mathbb{R}_{\geq 0}$. The parameters $b_1$, $b_2$ are selected such as to (safely) over-approximate $d_{\text{safe},ij}$ by over-approximating the nominator and under-approximating the denominator in \eqref{eq:CBF_CA_formCBF_dSafe}, see \prettyref{fig:problem_funcapprox} for further details. By substituting the smooth approximation \eqref{eq:CBF_CA_formCBF_maxApprox} in \eqref{eq:CBF_CA_formCBF_axMinEff}--\eqref{eq:CBF_CA_formCBF_axMinRelative}, we obtain the continuously differentiable safety distance $\tilde{d}_{\text{safe},ij}$ along with the collision avoidance CBF
\begin{align}
{h}_{c,ij}(x) \defeq d_{ij} - \tilde{d}_{\text{safe},ij}
\label{eq:CBF_CA_formCBF_hCA}
\end{align}
and the associated CBF constraint 
\begin{align}
L_{f} {h}_{c,ij}(x) + L_{g} {h}_{c,ij}(x) u  + \lambda_{c} {h}_{c,ij}(x) &\geq 0 \label{eq:CBF_CA_formCBF_CBFCA} 
\end{align}
where $\lambda_{c} \in \mathbb{R}_{>0}$ is a positive constant. 
\begin{proposition}
	The continuous differentiable function $h_{c,ij}(x)$ in \eqref{eq:CBF_CA_formCBF_hCA} is a CBF on $\mathcal{C}_{c,ij}$ for the dynamical system \eqref{eq:problem_modeling_lumpedSS}. Thus, the set $K_{c,ij}(x)$ of controls s.t. $u\in \mathcal{U}$ that render the zero superlevel set $\mathcal{C}_{c,ij}$ forward invariant is non-empty.
	\label{prop:CBFCA_formCBF_propIsCBF}
\end{proposition}
\vspace*{-5mm}
\begin{proof}
	By construction of CBF $h_{c,ij}(x)$, there exists a $u\in\mathcal{U}$ that renders  \eqref{eq:CBF_CA_formCBF_CBFCA}
	feasible for all $x\in\mathcal{C}_{c,ij}$. As such $h_{c,ij}(x)$ is a CBF on $\mathcal{C}_{c,ij}$. \hfill$\square$
\end{proof}
\begin{figure}[b!]		
	\centering
	\vspace*{-4mm}
	\begin{tikzpicture}[ auto, shorten >=2pt, shorten <=1pt, node distance=1.5cm, every loop/.style={},
thick,main node/.style={circle,draw=mydarkblue,font=\sffamily\bfseries,fill=myblue!10!}]

\node[main node] (1) {\color{myblue}2};
\node[main node] (2) [below left of=1] {\color{myred}1};
\node[main node] (3) [below right of=2] {\color{mycyan}4};
\node[main node] (4) [below right of=1] {\color{mygreen}3};

\path[every node/.style={font=\sffamily\small},color=mydarkblue,ultra thick]
(1) edge node [left] {} (4)
(2) edge node [right] {} (1)
(3) edge node [right] {} (2)
(4) edge node [left] {} (3);
\end{tikzpicture}

%
	\vspace*{-2mm}
	\caption{Undirected conflict graph $\mathcal{G}=(\mathcal{V},\mathcal{E})$ with $\mathcal{V} = \mathcal{A}$. 
	An edge $(i,j) \in \mathcal{E}$ indicates that Agent $i$ and Agent $j$ may collide along their paths.}
	\label{fig:problem_conflictGraph}
\end{figure}
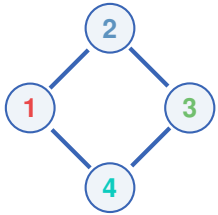

\vspace*{1mm}
In the centralized control scheme, it suffices to either impose CBF $h_{c,ij}(x)$ or $h_{c,ji}(x)$ for two conflicting agents $i,j \in \mathcal{A}$. To this end, we define the conflict graph $\mathcal{G}=(\mathcal{V},\mathcal{E})$ with vertices $\mathcal{V} = \mathcal{A}$ and edges $\mathcal{E} \subseteq \mathcal{V} \times \mathcal{V}$. An edge \mbox{$(i,j) \in \mathcal{E}$} with $i\neq j$ indicates that Agent $i$ and \mbox{Agent $j$} may collide along their paths. An example conflict graph for the scenario in 	\prettyref{fig:problem_description_intersectionScenario} is depicted in \prettyref{fig:problem_conflictGraph}. We finally define the conflict set of Agent $i$ as 
\(
\mathcal{A}_c(i) \defeq \{ j\in\mathcal{A} \mid (i,j)\in\mathcal{E} \land j>i \}
\) and impose \eqref{eq:CBF_CA_formCBF_CBFCA} for all $j\in\mathcal{A}_c(i)$.
\begin{remark}
	CBF \eqref{eq:CBF_CA_formCBF_hCA} is tailored to the straight crossing scenario but can be further extended to accommodate additional scenarios, such as those that involve turning agents. 
\end{remark}


\subsection{Control-sharing Property}
\label{sec:CBFCA_csProperty}

\begin{proposition}
	For every agent $i$ and every two conflicting agents $j,l \in \mathcal{A}_c(i)$, the CBFs $h_{c,ij}(x)$ and $h_{c,il}(x)$ exhibit a control-sharing property, i.e., $K_{{c,ij}}(x) \cap K_{{c,il}}(x) \neq \emptyset$. 
\end{proposition}	
\vspace*{-1mm}
\begin{proof}
  Following \prettyref{sec:CBF_backgrnd}, we start within the safe set, i.e., $x(t_0) \in \mathcal{C}_{c,ij}$ and $x(t_0) \in \mathcal{C}_{c,il}$. 
  By construction of \eqref{eq:CBF_CA_formCBF_hCA} we know that there exist control actions $(u_i,u_j)$ such that a collision between Agent $i$ and \mbox{Agent $j\in\mathcal{A}_c(i)$} can be avoided. 
  Likewise, there exist control actions $(u_i^\prime,u_l)$ to avoid a collision between Agent $i$ and Agent $l\in\mathcal{A}_c(i)$. In both cases, (at least) by choosing $u_i=u_i^\prime=\underline{a}_{x,\text{eff},i}$ (see \eqref{eq:CBF_CA_formCBF_axMinEff}) collision avoidance can be guaranteed (of course, there may be other choices too) --- then, $u_j \in \mathcal{U}_j$ and $u_l \in \mathcal{U}_l$ can be selected such as to satisfy CBF constraint \eqref{eq:CBF_CA_formCBF_CBFCA}. Consequently, it holds $K_{{c,ij}}(x) \cap K_{{c,il}}(x) \neq \emptyset$. \hfill$\square$
\end{proof}

\begin{proposition} \label{prop:CBFCA_csProperty_propCAandVel}
	Collision avoidance CBFs \eqref{eq:CBF_CA_formCBF_hCA} 
	exhibit a control-sharing property with velocity CBFs \eqref{eq:CBF_velCBF_hLowerBnd} and \eqref{eq:CBF_velCBF_hUpperBnd}. In other words, for Agent $i$ and Agent $j \in \mathcal{A}_c(i)$ \eqref{eq:CBF_CA_formCBF_hCA} is jointly feasible with \eqref{eq:CBF_velCBF_hLowerBnd}, \eqref{eq:CBF_velCBF_hUpperBnd} when utilizing $h_{\underline{v},i}(x_i)$, $h_{\overline{v},i}(x_i)$ for Agent $i$ and $h_{\underline{v},j}(x_j)$, $h_{\overline{v},j}(x_j)$ for Agent $j$.	
\end{proposition}
\vspace*{-1mm}
\begin{proof} 
	By virtue of \prettyref{sec:CBF_backgrnd}, we start within the safe set, i.e., 
	$x_i(t_0) \in \mathcal{C}_{\underline{v},i}$, 
	$x_i(t_0) \in \mathcal{C}_{\overline{v},i}$,  
	$x_j(t_0) \in \mathcal{C}_{\underline{v},j}$, 
	$x_j(t_0) \in \mathcal{C}_{\overline{v},j}$ and 
	$x(t_0) \in \mathcal{C}_{c,ij}$. 
	Let's assume that there exist initial states $x_i(t_0)$, $x_j(t_0)$ such that velocity CBF constraints \eqref{eq:CBF_velCBF_CBFconsLowerBnd} and \eqref{eq:CBF_velCBF_CBFconsUpperBnd} are satisfied while collision avoidance CBF constraint \eqref{eq:CBF_CA_formCBF_CBFCA} is not. 
	By construction of CBF \eqref{eq:CBF_CA_formCBF_hCA}, a collision can always be avoided when both agents apply their maximum effective deceleration
	\mbox{$\underline{a}_{x,\text{eff},i}=\max\{\underline{a}_{x,i},\,-\lambda_{\underline{v}} v_i \}$} 
	and \mbox{$\underline{a}_{x,\text{eff},j}=\max\{\underline{a}_{x,j},\,-\lambda_{\underline{v}} v_j \}$}. 
	By virtue of \eqref{eq:CBF_CA_formCBF_axMinEff}, the maximum effective deceleration of both agents is compatible with the lower acceleration bound \eqref{eq:CBF_velCBF_vminLowerAccel} induced by velocity CBF \eqref{eq:CBF_velCBF_hLowerBnd}. 	
	By construction, the upper velocity bound is never in conflict with collision avoidance CBF \eqref{eq:CBF_CA_formCBF_hCA}.
	That said, at all times there exists at least one control action $(u_i,u_j)$ for \mbox{Agent $i$} and Agent $j$ that is jointly feasible with CBFs \eqref{eq:CBF_velCBF_hLowerBnd}, \eqref{eq:CBF_velCBF_hUpperBnd} and \eqref{eq:CBF_CA_formCBF_hCA}. This statement, though, is in contradiction with the assumption that \eqref{eq:CBF_CA_formCBF_CBFCA} is jointly infeasible.  \hfill$\square$
\end{proof}

\section{Controller Synthesis}
\label{sec:synthesis}

To synthesize the control law for the nonlinear system \eqref{eq:problem_modeling_xdot}, the State-Dependent Riccati Equation (SDRE) method \cite{Cloutier1997a} is applied. By virtue of \prettyref{prob:problem_description_problemStatement}, the control objective of every agent $i$ is to track its reference speed $v_{\text{ref},i}$. As such, for control purposes $s_i$ can be removed from the state space while adding the state equation $\dot{e}_i = v_{\text{ref},i}-v_i$ is required to incorporate integral action. That said, with the state vector $\xi_i=[v_i,~e_i]^{\T}$ and control input $u_i=a_{x,i}$ we can rewrite the nonlinear dynamics \eqref{eq:problem_modeling_xdot} of Agent $i$ as
\vspace*{-1mm}
\begin{align}
\dot{\xi}_i = A_i(\xi_i) \xi_i + B_i u_i + E_i v_{\text{ref},i} ~~~~~~~~~~~~
\label{eq:synthesis_ssSDRE} \raisetag{6mm}\\
\text{with}~~
A_i(\xi_i) \defeq \begin{bmatrix}
-a_{11,i} & 0 \\
-1 & 0
\end{bmatrix},~
B_i \defeq \begin{bmatrix}
1 \\ 0
\end{bmatrix},~
E_i \defeq \begin{bmatrix}
0 \\ 1
\end{bmatrix} \notag
\end{align}
where the system matrix $A_i(\xi_i)$ depends on the state $\xi_i$ while 
$a_{11,i}=F_{r,i}(v_i)/(m_i v_i)$ if $v_i \geq v_{\text{thld}}$ and $a_{11,i}=0$ otherwise. Moreover, $v_{\text{thld}}=0.1$ is a lower velocity threshold to avoid division by zero.  
For every agent $i$, the control action 
is obtained by minimizing 
\begin{align}
J = \frac{1}{2}\,\int_{0}^{\infty}  \xi_{i}^{\T} Q_i \,\xi_{i} + u_{i}^{\T} R_i\, u_{i}\, dt \label{eq:synthesis_costDSDRE}\\[-7mm] \notag
\end{align}
with $Q \succeq 0$ and $R \succ 0$, given that $(A_{i}(\xi_i),B_{i})$ is pointwise controllable for every measurable $\xi_{i}$ ($v_i$ and $e_i$ are considered to be measurable). The solution of the infinite horizon optimal control problem can be obtained by solving the SDRE \cite{Cloutier1997a} in dependence of the current state $\xi_i$. In our case, we can even pre-compute an analytical expression for the solution $P_i(\xi_i)$ of the SDRE. Then, the feedback gain $K_i$ can efficiently be calculated by $K_i = R_i^{-1}B_i^{\T}P_i(\xi_i)$ and the \textit{nominal} control action is given as
\begin{align}
u_{\text{nom},i} &= -K_{i} \begin{bmatrix}
	v_{i} - v_{\text{ref},i} \\
	e_i
\end{bmatrix}. \label{eq:synthesis_discControlLaw} \\[-6.5mm] \notag
\end{align}
We utilize the CBFs in \prettyref{sec:CBF} and \prettyref{sec:CBFCA} to complement 
control law \eqref{eq:synthesis_discControlLaw} with constraints. 
At time $t$ these CBFs prescribe a set of feasible controls that render the system safe at all times $t^\prime\geq t$. 
That said, given the nominal control action $u_{\text{nom}}(t) \defeq \{u_{\text{nom},i}(t)\}_{i=1}^{N_\mathcal{A}}$ and measured state $x(t)$ (see Assumption A4) 
at time $t$, we solve the following convex QP problem which minimizes the 2-norm of the difference 
between the corrected control action $u$ and $u_{\text{nom}}$
subject to linear CBF and input constraints $u\in\mathcal{U}$, i.e.,

\vspace*{2mm}
\noindent
\textbf{Centralized CBF-QP}
\vspace*{-1mm}
\begin{subequations}\label{eq:synthesis_minNormQP}
	\begin{align} 
	\hspace*{-2mm}u(x) =~& {\mathargmin_{u \in \mathcal{U}}}~ \frac{1}{2} \| u - u_{\text{nom}} \|_2^2 \\
	\hspace*{-2mm}\mathst~~~    &~\hspace*{-3mm} \textbf{velocity CBFs:}~~\forall i \in \mathcal{A} \notag\\
	\hspace*{-2mm} &~\hspace*{-3mm}L_{f_i} {h}_{\underline{v},i}(x_i) + L_{g_i} {h}_{\underline{v},i}(x_i) u_i  + \lambda_{\underline{v}} {h}_{\underline{v},i}(x_i) \geq 0 \raisetag{4.2mm}\\
	\hspace*{-2mm}&~\hspace*{-3mm}L_{f_i} {h}_{\overline{v},i}(x_i) + L_{g_i} {h}_{\overline{v},i}(x_i) u_i  + \lambda_{\overline{v}} {h}_{\overline{v},i}(x_i) \geq 0\\[1mm]			
	\hspace*{-2mm} &~\hspace*{-3mm} \textbf{collision avoidance CBFs:~~}\forall i \in \mathcal{A},~ \forall j \in \mathcal{A}_c(i)\notag\\
	\hspace*{-2mm}&~\hspace*{-3mm}L_f {h}_{c,ij}(x) + L_g {h}_{c,ij}(x) u  + \lambda_{c} h_{c,ij}(x) \geq 0
	\end{align}
\end{subequations}

\begin{theorem}
	The centralized CBF-QP \eqref{eq:synthesis_minNormQP} solves \prettyref{prob:problem_description_problemStatement} and renders the system provably safe. Moreover, the CBFs are jointly feasible under input constraints. 
\end{theorem}
\vspace*{-1mm}
\begin{proof}
	Follows from \prettyref{prop:CBF_velCBF_propCBFvel} -- \prettyref{prop:CBFCA_csProperty_propCAandVel}. \hfill$\square$
\end{proof}

\section{Simulation Results}
\label{sec:sim}

\begin{figure}[b!]
	\vspace*{-6mm}
	\setlength\fwidth{0.42\textwidth}	
	\input{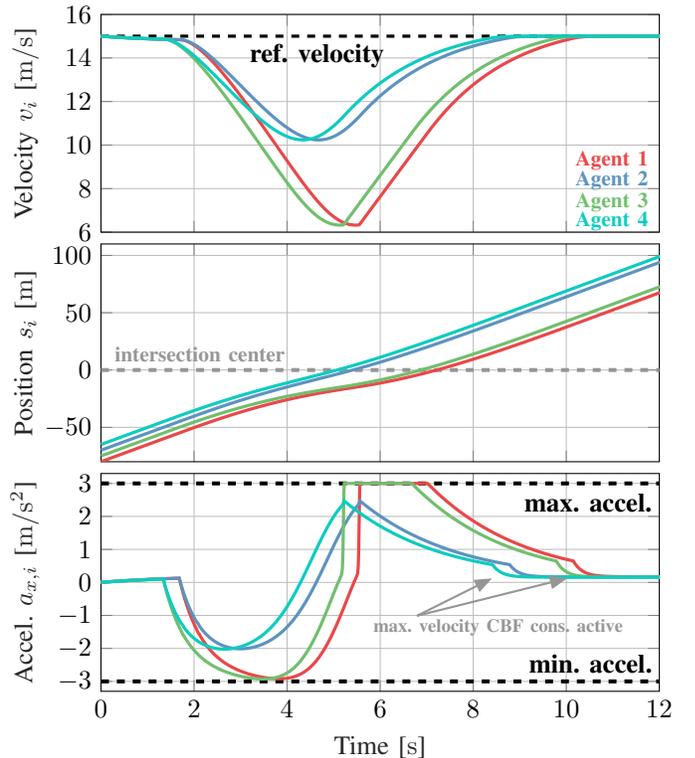}	
	\vspace*{-8mm}
	\caption{States (velocity and path position) and inputs (acceleration) of Agent 1 ({\color{myred}red}), Agent 2 ({\color{myblue}blue}), Agent 3 ({\color{mygreen}green}) and Agent 4 ({\color{mycyan}cyan}).} 
	\vspace*{-3mm}
	\label{fig:sim_results_statesInputs}
\end{figure}
\begin{figure*}
	\centering
	\setlength\fwidth{0.70\textwidth}	
	\input{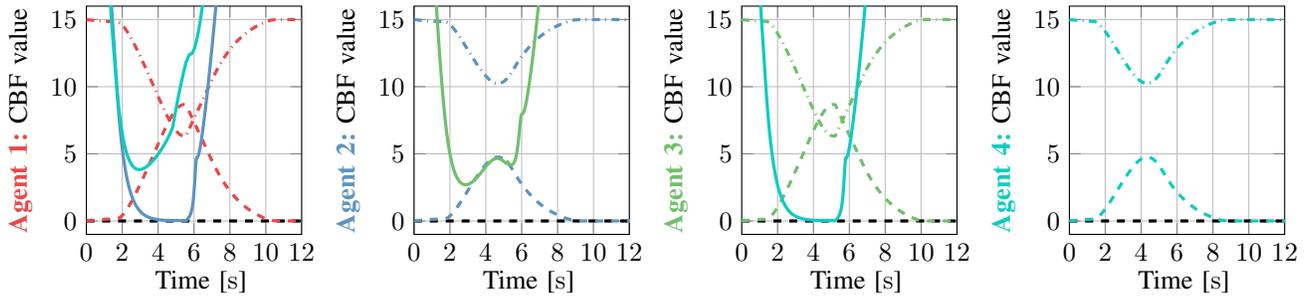}
	\vspace*{-4mm}
	\caption{CBF values for Agent $i$ in column $i$. The dashed and dotted-dashed lines in the color of Agent $i$ correspond to the maximum and minimum velocity CBFs, respectively. Solid lines indicate the collision avoidance CBF values in the color of the conflicting agent.}		
	\vspace*{-6mm}
	\label{fig:sim_results_CBF}
\end{figure*}

\subsection{Simulation Setup}
\label{sec:sim_setup}

We investigate 
a realistic intersection scenario with four straight crossing agents as depicted in \prettyref{fig:problem_description_intersectionScenario}. Their conflict graph corresponds to \prettyref{fig:problem_conflictGraph}. All agents have dimensions $L_i=\unit[5]{m}$ and $W_i=\unit[2]{m}$, and their mass is: $m_1=\unit[1200]{kg}$, $m_2=\unit[1300]{kg}$, $m_3=\unit[1400]{kg}$, $m_4=\unit[1500]{kg}$. 
For the different agent masses, the resistance force coefficients are selected as: $c_{0,i}=0.01 m_i g$ (with gravitational acceleration $g$), $c_{1,i}=-0.433$ and \mbox{$c_{2,i}=0.422$} \cite{Kuehlwein2016a}. 
The initial and reference velocities for all agents are $\unitfrac[15]{m}{s}$. Moreover, their initial positions in the global reference frame are: 
$(-80,-2)$ for Agent 1, 
$(-2,70)$ for Agent 2, 
$(75,2)$ for Agent 3 and 
$(2,-65)$ for Agent 4.
For $s_i=0$, the agents are at the center of the intersection such that either $X_i$ or $Y_i$ are zero, see \prettyref{fig:problem_description_intersectionScenario}. The weights in SDRE \eqref{eq:synthesis_costDSDRE} are the same for all agents, i.e., \mbox{$Q=\mathrm{diag}(1,0.05)$} and $R=4$.  For the CBFs, we have chosen the following constants: $\lambda_{c}=2$, $\lambda_{\underline{v}}=5$, $\lambda_{\overline{v}}=5$. The superellipse parameters $\Delta a_{i}$ and $\Delta b_{i}$ are set to $\unit[1.5]{m}$ (safety buffer) for all agents. 
Moreover, the longitudinal acceleration of every agent is constrained by 
$u_i \in [-3, 3] \,\unitfrac{m}{s^2}$ 
and the maximum velocity $\overline{v}_i$ is set to $\unitfrac[15]{m}{s}$ $(=v_{\text{ref},i})$. Simulations are run on an Intel i7 machine at $\unit[2.5]{GHz}$ with Matlab R2021b. CBF-QP \eqref{eq:synthesis_minNormQP} is solved every 
$\unit[10]{ms}$ 
using CasADi v3.5.5 with qpOASES as QP solver \cite{Andersson2019a}. 
To avoid model uncertainties (see  \prettyref{ass:problem_description_problemAssump}), the simulation and control-oriented model are the same. 

\subsection{Discussion of Results}
\label{sec:sim_discussion}

\prettyref{fig:sim_results_statesInputs} conveys the resulting states (top: velocity, middle: path position) and inputs (bottom: acceleration) for the given scenario. It can be recognized that all agents reduce their speed when approaching the intersection, which is in line with choosing the same weights $Q,R$ for all agents. 
\mbox{Agent 2} (blue) and Agent 4 (cyan) cross the intersection first with a velocity of $\unitfrac[10.2]{m}{s}$ before Agent 1 (red) and Agent 3 (green) are able to cross too. 
When decelerating for agents 2 and 4, Agent 1 and Agent 3 reach their lower acceleration limits and a minimum velocity of $\unitfrac[6.3]{m}{s}$. Due to the slightly conservative $\max$-operator approximation in \prettyref{sec:CBF_CA_formCBF}, the acceleration stays marginally above the minimum. During that time, the hard constrained CBF-QP \eqref{eq:synthesis_minNormQP} remains feasible and input constraints are not violated. After crossing the intersection, every agent accelerates to track its reference velocity. When being very close to the set-point, the maximum velocity CBF constraint  becomes active and reduces the acceleration smoothly to avoid overshooting. 

\prettyref{fig:sim_results_CBF} provides further insights into the function values of the utilized CBFs for every agent. The dashed (\tikz[baseline=-0.5ex]\draw [thick,dashed] (0,0) -- (0.3,0);) and dotted-dashed (\tikz[baseline=-0.5ex]\draw [thick,dash dot] (0,0) -- (0.4,0);) lines in the color of Agent $i$ in column $i$ (e.g. red dashed and dotted-dashed for \mbox{Agent 1} in \mbox{column 1)} show the CBF values for the maximum and minimum velocity CBF, respectively. The solid lines represent the collision avoidance CBFs in the color of the conflicting Agent $j \in \mathcal{A}_c(i)$ (e.g. the cyan solid line in column 1 is CBF $h_{c,14}$ to avoid collisions between agents 1 and 4). Evidently, every CBF is greater than or equal to zero (black dashed line), i.e., the system is safe at all times. 

The centralized CBF-QP \eqref{eq:synthesis_minNormQP} has been solved in real-time, i.e., within $\unit[0.7]{ms}$ on average and $\unit[3.0]{ms}$ in the worst case. 
Concluding, the proposed control scheme solves \prettyref{prob:problem_description_problemStatement} and renders the system safe subject to input constraints. 

\section{Conclusions}
\label{sec:conclusion}

We have proposed a provably safe centralized CBF-based control framework to coordinate connected agents at intersections. The utilized CBFs exhibit a control-sharing property while being compatible with input constraints. A simulation study has demonstrated the efficacy and computational efficiency of our concept.
A limitation of CBFs, though, is their limited ability to exploit predictions (e.g. when other non-controllable road users come into play). Moreover, CBFs may be prone to deadlock-like situations, e.g., due to geometric symmetry in the initial condition. 
 Therefore, we envision to investigate a combination of MPC- and CBF-based control schemes. Additionally, we intend to study a  distributed numerical solution of CBF-QP \eqref{eq:synthesis_minNormQP}.




\bibliographystyle{IEEEtran}
\bibliography{IEEEabrv}

\begin{thebibliography}{10}
\providecommand{\url}[1]{#1}
\csname url@samestyle\endcsname
\providecommand{\newblock}{\relax}
\providecommand{\bibinfo}[2]{#2}
\providecommand{\BIBentrySTDinterwordspacing}{\spaceskip=0pt\relax}
\providecommand{\BIBentryALTinterwordstretchfactor}{4}
\providecommand{\BIBentryALTinterwordspacing}{\spaceskip=\fontdimen2\font plus
\BIBentryALTinterwordstretchfactor\fontdimen3\font minus
  \fontdimen4\font\relax}
\providecommand{\BIBforeignlanguage}[2]{{%
\expandafter\ifx\csname l@#1\endcsname\relax
\typeout{** WARNING: IEEEtran.bst: No hyphenation pattern has been}%
\typeout{** loaded for the language `#1'. Using the pattern for}%
\typeout{** the default language instead.}%
\else
\language=\csname l@#1\endcsname
\fi
#2}}
\providecommand{\BIBdecl}{\relax}
\BIBdecl

\bibitem{Khayatian2020a}
M.~Khayatian~et al., ``{A Survey on Intersection Management of Connected
  Autonomous Vehicles},'' \emph{ACM Transactions on Cyber-Physical Systems},
  vol.~4, no.~4, pp. 1--27, 2020.

\bibitem{Malikopoulos2021a}
A.~A. Malikopoulos, L.~Beaver, and I.~V. Chremos, ``{Optimal time trajectory
  and coordination for connected and automated vehicles},'' \emph{Automatica},
  vol. 125, p. 109469, 2021.

\bibitem{Katriniok2019a}
A.~Katriniok, P.~Sopasakis, M.~Schuurmans, and P.~Patrinos, ``{Nonlinear Model
  Predictive Control for Distributed Motion Planning in Road Intersections
  Using PANOC},'' in \emph{IEEE Conference on Decision and Control}, 2019, pp.
  5272--5278.

\bibitem{Ames2016a}
A.~D. Ames, X.~Xu, J.~W. Grizzle, and P.~Tabuada, ``{Control Barrier Function
  Based Quadratic Programs for Safety Critical Systems},'' \emph{Trans. on
  Automatic Control}, vol.~62, no.~8, pp. 3861--3876, 2017.

\bibitem{Ames2019a}
A.~D. Ames, S.~Coogan, M.~Egerstedt, G.~Notomista, K.~Sreenath, and P.~Tabuada,
  ``{Control Barrier Functions: Theory and Applications},'' in \emph{European
  Control Conference}, 2019, pp. 3420--3431.

\bibitem{Katriniok2017a}
A.~Katriniok, P.~Kleibaum, and M.~Jo\v{s}evski, ``{Distributed Model Predictive
  Control for Intersection Automation Using a Parallelized Optimization
  Approach},'' in \emph{IFAC World Congress}, vol.~50, no.~1, 2017, pp.
  5940--5946.

\bibitem{Wang2017a}
L.~Wang, A.~D. Ames, and M.~Egerstedt, ``{Safety Barrier Certificates for
  Collisions-Free Multirobot Systems},'' \emph{IEEE Transactions on Robotics},
  vol.~33, no.~3, pp. 661--674, 2017.

\bibitem{Squires2018a}
E.~Squires, P.~Pierpaoli, and M.~Egerstedt, ``{Constructive Barrier
  Certificates with Applications to Fixed-Wing Aircraft Collision Avoidance},''
  in \emph{IEEE Conf. on Control Techn. \& Appl.}, 2018, pp. 1656--1661.

\bibitem{Xiao2019a}
W.~Xiao, C.~Belta, and C.~G. Cassandras, ``Decentralized merging control in
  traffic networks: a control barrier function approach,'' in \emph{Conference
  on Cyber-Physical Systems}, 2019, pp. 270--279.

\bibitem{Santillo2021a}
M.~Santillo and M.~Jankovic, ``{Collision Free Navigation with Interacting,
  Non-Communicating Obstacles},'' in \emph{IEEE American Control Conference},
  2021, pp. 1637--1643.

\bibitem{Khaled2020a}
S.~Khaled, O.~M. Shehata, and E.~I. Morgan, ``{Intersection Control for
  Autonomous Vehicles Using Control Barrier Function Approach},'' in
  \emph{Novel Intell. \& Leading Emerging Sciences Conf.}, 2020, pp. 479--485.

\bibitem{Shivam2020a}
S.~Shivam, Y.~Wardi, M.~Egerstedt, A.~Kanellopoulos, and K.~Vamvoudakis,
  ``{Intersection-Traffic Control of Autonomous Vehicles using Newton-Raphson
  Flows and Barrier Functions},'' in \emph{IFAC World Congress}, vol.~53,
  no.~2, 2020, pp. 15\,733--15\,738.

\bibitem{Agrawal2021a}
D.~Agrawal and D.~Panagou, ``{Safe Control Synthesis via Input Constrained
  Control Barrier Functions},'' in \emph{IEEE Conference on Decision and
  Control}, 2021, pp. 6113--6118.

\bibitem{Xu2018a}
X.~Xu, ``{Constrained control of input-output linearizable systems using
  control sharing barrier functions},'' \emph{Automatica}, vol.~87, pp.
  195--201, 2018.

\bibitem{Rajamani2012}
R.~Rajamani, \emph{{Vehicle Dynamics and Control}}.\hskip 1em plus 0.5em minus
  0.4em\relax Springer, 2012, vol.~2.

\bibitem{Khalil2002a}
H.~K. Khalil, \emph{{Nonlinear Systems}}, 3rd~ed.\hskip 1em plus 0.5em minus
  0.4em\relax Prentice Hall, 2002.

\bibitem{Breeden2021a}
J.~Breeden and D.~Panagou, ``{High Relative Degree Control Barrier Functions
  Under Input Constraints},'' in \emph{IEEE Conference on Decision and
  Control}, 2021, pp. 6119--6124.

\bibitem{Cloutier1997a}
J.~Cloutier, ``State-dependent riccati equation techniques: an overview,'' in
  \emph{American Control Conference}, 1997, pp. 932--936.

\bibitem{Kuehlwein2016a}
J.~K\"{u}hlwein, ``{Driving resistances of light-duty vehicles in Europe:
  Present situation, trends, and scenarios for 2025},'' in \emph{The
  International Council on Clean Transportation (White Paper)}, 2016.

\bibitem{Andersson2019a}
J.~A.~E. Andersson, J.~Gillis, G.~Horn, J.~B. Rawlings, and M.~Diehl,
  ``{CasADi} -- {A} software framework for nonlinear optimization and optimal
  control,'' \emph{Math. Programming Computation}, pp. 1--36, 2019.

\end{thebibliography}


%

%


\end{document}